\documentclass[english,a4paper, 12pt,dvipdfm]{amsart}
\makeindex
\usepackage[all,web]{xy}
\usepackage{graphicx}
\usepackage{amsmath} 
\usepackage{amsfonts}
\usepackage{amssymb}
\usepackage{amsthm}
\usepackage{setspace}

 \theoremstyle{plain}    
 \newtheorem{thm}{Theorem}[section]
 \numberwithin{equation}{section} 
 \numberwithin{figure}{section} 
 \theoremstyle{remark}   
  
 \theoremstyle{plain}    
 \theoremstyle{plain}    
  
 \theoremstyle{remark} 
 
 \theoremstyle{remark}
 
 \theoremstyle{definition}
  
 \theoremstyle{plain}
 \newtheorem{conjecture}[thm]{Conjecture} 
 \theoremstyle{plain}    
 \theoremstyle{plain}    
 \theoremstyle{definition}
 
 \theoremstyle{definition}
  \newtheorem{example}[thm]{Example}
 \theoremstyle{plain}    
  
 \theoremstyle{plain}    
 \newtheorem{lem}[thm]{Lemma} 
 \theoremstyle{remark}    
  
 \theoremstyle{remark}    
  
 \theoremstyle{definition}
  
 \theoremstyle{plain}
 \theoremstyle{remark}
 
 \theoremstyle{remark}

 1

\def\PP{\mathbb{P}}
\def\QQ{\mathbb{Q}}

\def\ZZ{\mathbb{Z}}
\def\hyper{\mathcal H}

\newcommand{\vtwo}[2]{\left(
        \begin{matrix}#1\\#2
        \end{matrix}\right)}
\newcommand{\mtwo}[4]{\left(
        \begin{matrix}#1&#2\\#3&#4
        \end{matrix}\right)}
\newcommand{\abcd}[4]{\left(
        \begin{smallmatrix}#1&#2\\#3&#4\end{smallmatrix}\right)}

\def\inner<#1,#2>{{\left\langle{{#1},{#2}}\right\rangle}}
\def\sl2of#1{\textrm{SL}_2(#1)}

\def\inner<#1,#2>{{\left\langle{{#1},{#2}}\right\rangle}}

\usepackage{babel}

\begin{document}

\title{Further properties of a function of Ogg and Ligozat}
\author{Carlos Casta\~no-Bernard}

\begin{abstract}
Certain identities of Ramanujan may be succinctly expressed in terms of the
rational function $\breve{g}_\chi=\breve{f}_\chi-\frac{1}{\breve{f}_\chi}$
on the modular curve $X_0(N)$,
where $\breve{f}_\chi=w_Nf_\chi$ and $f_\chi$ is a certain modular unit
on the Nebentypus cover $X_\chi(N)$ introduced by Ogg and Ligozat
for prime $N\equiv 1\pmod{4}$ and $w_N$ is the Fricke involution.
These correspond to levels $N=5,13$,
where the genus $g_N$ of $X_0(N)$ is zero.
In this paper we study a slightly more general kind of relations
for each $\breve{g}_\chi$ such that $X_0(N)$ has genus $g_N=1,2$,
and also for each $h_\chi=g_\chi + \breve{g}_\chi$
such that the Atkin-Lehner quotient $X_0^+(N)$ has genus $g_N^+=1,2$.
It turns out that if $n$ is the degree of 
the field of definition $F$ of the non-trivial zeros of the latter,
then the degree of the normal closure of $F$ over $\QQ$ 
is the $n$-th solution of Singmaster's Problem.
\end{abstract}


\maketitle

\pagenumbering{roman}
\setcounter{page}{0}



\tableofcontents


\pagenumbering{arabic}
\pagestyle{headings}

%
%
%
\section{Introduction}
Let $X_0(N)$ be the usual compactification of the
coarse moduli space $Y_0(N)$ of
isomorphism classes of pairs $(E, E^\prime)$
of elliptic curves $E$ and $E^\prime$
linked by a cyclic isogeny of degree $N$.
Suppose $N$ is a prime number such that $N\equiv 1\pmod{4}$
and let $X_\chi(N)\longrightarrow X_0(N)$ be the Nebentypus covering
that corresponds to the kernel of the
quadratic character $\chi$ on $(\ZZ/N\ZZ)^\times$
as in (7.5.14) and (7.5.15) of Shimura's book~\cite{shimura:automorphic}.
In Mazur's paper~\cite[pp.~107,~108]{mazur:eisen}
it may be found a construction due to Ligozat of
a modular unit $f_\chi$ on $X_\chi(N)$.
Certain identities due to Ramanujan may be succinctly expressed 
in terms of the rational
function $\breve{g}_\chi = \breve{f}_\chi-\frac{1}{\breve{f}_\chi}$
in $\QQ(X_0(N))$ as $\breve{g}_\chi(\tau)=P(t(\tau))$, 
where $t(\tau)$ is a certain modular unit on $X_0(N)$, $P(T)\ZZ[T]$
is a polynomial of degree $1$, $\breve{f}_\chi=w_Nf_\chi$,
and $w_N$ is the Fricke involution.
(See Section~\ref{sec:ram}.)
In these cases the genus $g_N$ of $X_0(N)$ is $g_N=0$.
It turns out that for $X_0(N)$ of genus $g_N>0$,
identities of the form $\breve{g}_\chi = P(t(\tau))$
with $P(T)\QQ[T]$ arbitrary are unlikely.
So the natural analogs of these identities for genus $g_N>0$ 
are perhaps of the form $\breve{g}_\chi = P(X,Y)$,
for some ``canonically defined'' rational functions $X$ and $Y$ of $X_0(N)$.
We suggest such identities for $\breve{g}_\chi$ for $g_N=1$ and $2$,
and also for the
rational function $h_\chi = g_\chi + \breve{g}_\chi\in\QQ(X_0^+(N))$,
where $g_\chi=f_\chi - \frac{1}{f_\chi}$
and $X_0^+(N) = X_0(N)/\{1,w_N\}$
is the Atkin-Lehner quotient of $X_0(N)$ defined by the involution $w_N$
and $X_0^+(N)$ has genus $g_N^+=1$ and $2$.
With the help of some Gr\"obner basis algorithms,
the latter identities yield
the field of definition of each of the zeros of $h_\chi$.
These turn out to be either $\QQ$ or
a finite extension $F$ of the real quadratic field $\QQ(\sqrt{N})$;
the field extension $F/\QQ$ has degree $n=\frac{1}{2}B_{2,\chi}$,
if $N\equiv 1\pmod{8}$ and $n=\frac{1}{2}B_{2,\chi}-1$,
if $N\equiv 5\pmod{8}$,
where $B_{m,\chi}$ is the $m$-th generalised Bernoulli number
attached to the character $\chi$.
Moreover,
in each of these cases the Galois group $G=G(F^\textit{nrm}/\QQ)$ of
the normal closure $F^\textit{nrm}$ of $F$ over the rationals $\QQ$ is
the wreath product $G=S_{\frac{n}{2}}\wr C_2$
of the symmetric group $S_{\frac{n}{2}}$
of permutations of $\frac{n}{2}$ objects
and the cyclic group of order two $C_2$.
This means that in each of these cases
the Galois group $G$  has order $\#(G)=2(\frac{n}{2})!$,
which turns out to be 
the $n$-th solution of Singmaster's Problem~\cite{singmaster:problem}.

\bigskip
\textbf{Organisation of the paper.}
In order to make the exposition as self contained as possible
we include some standard results on modular curves and modular units
tailored to our needs,
in Section~\ref{sec:back}.
We also include a variant of Ligozat's construction based on
an extension of a classical identity used by Gau\ss\ in his third proof
of the Law of Quadratic Reciprocity,
in Section~\ref{sec:lig}.
The main results are contained in
Section~\ref{sec:ram},
the $g_N=1$ case (i.e. $N = 17$)
and the $g_N=2$ case (i.e. $N = 29, 37$),
and in Section~\ref{sec:hchi},
the field of definition of the zeros of $h_\chi(\tau)$
for $g_N^+=1$ (i.e. $N = 37, 53, 61, 89, 101$),
and for $g_N^+=2$ (i.e. $N = 73$).
The paper concludes with a conjecture and an open problem,
in Section~\ref{sec:conj}.

\bigskip
\textbf{Acknowledgements.}
Our computer calculations were performed with the help of the
computer algebra packages
\textsc{Magma}~\cite{bosma:magma},
\textsc{Macaulay2}~\cite{grayson:m2},
and \textsc{Pari}~\cite{pari:gp}
installed on GNU/Linux computers.
We wish to thank
the Department of Pure Mathematics and Mathematical Statistics
of the University of Cambridge for the access to these packages.
We also wish to express our gratitude to my wife,
Concepci\'on Vargas-Ch\'avez and my brother,
Ricardo Casta\~no-Bernard
for their constant support and encouragement.


%
\section{Background material}\label{sec:back}
Let $X_1(N)$ be the usual compactification of the coarse
moduli space $Y_1(N)$ of isomorphism classes of unordered pairs 
\begin{displaymath}
\{(E,P),(E,-P)\},
\end{displaymath}
where $E$ is an elliptic curve and $P$ a point of $E$.
Let $\pi$ be the natural degeneracy morphism
\begin{equation*}
Y_1(N)\longrightarrow Y_0(N)
\end{equation*}
induced by the map
\begin{equation*}
\{(E,P),(E,-P)\}\mapsto(E,E/\langle P\rangle),
\end{equation*}
where  $\langle P\rangle$ denotes the group generated by the point $P$.
Note that for each $r\in\ZZ$ such that $(r,N)=1$
the map $(E,p)\mapsto(E_{/F},r\cdot P)$
induces an automorphism $\sigma_r$ of $X_1(N)$ over $X_0(N)$.
Moreover,
the map $r\mapsto \sigma_r$ induces an isomorphism
from the multiplicative group $C_N=(\ZZ/N\ZZ)^\times/\{-1,1\}$
onto the Galois group $G(X_1(N)^\textit{an}/X_0(N)^\textit{an})$
of $X_1(N)^\textit{an}$ over $X_0(N)^\textit{an})$
Recall that
\begin{equation*}
X_0(N)^\textit{an}=\Gamma_0(N)\backslash\mathcal{H}^*,
\end{equation*}
and
\begin{equation*}
X_1(N)^\textit{an}=\Gamma_1(N)\backslash\mathcal{H}^*,
\end{equation*}
where
\begin{equation*}
\Gamma_0(N)=
\left\{\mu=\mtwo{\alpha}{\beta}{\gamma}{\delta}\in\textrm{SL}_2(\ZZ)
\,\colon\,
\gamma\equiv 0\pmod{N}\right\}.
\end{equation*}
and
\begin{equation*}
\Gamma_1(N)=
\left\{\mu=\mtwo{\alpha}{\beta}{\gamma}{\delta}\in\Gamma_0(N)\,\colon\,
\delta\equiv 1\pmod{N}\right\}.
\end{equation*}
The map
\begin{equation*}
\mtwo{\alpha}{\beta}{\gamma}{\delta}\mapsto \{\delta\}
\end{equation*}
induces the above isomorphism from
the set of right cosets $[\Gamma_0(N):\Gamma_1(N)]$
onto the multiplicative group $C_N$,
where
\begin{equation*}
\xymatrix{
\{\cdot\}\colon\ZZ \ar [r]& \left\{0,1,\dots,\frac{N-1}{2}\right\}\\
a \ar @{|->}[r]& \pm a\pmod{N},\\
}
\end{equation*}
following the notation of
Ogg~\cite{ogg:rational} and
Csirik's thesis~\cite{csirik:thesis}.
From now on we assume $N$ is a prime and such that $N\equiv 1\pmod{4}$,
so that $C_N$ is cyclic and
the kernel $\Omega=\textrm{ker}(\chi)$ of the quadratic character
\begin{equation*}
\xymatrix{
\chi\colon C_N\ar [r]& \{-1,1\},\\
n\ar @{|->}[r]& (\frac{n}{N}),\\
}
\end{equation*}
is a (cyclic) subgroup of index $2$.
The \textit{Nebentypus} curve $X_\chi(N)$ associated to
the character $\chi$ is
the intermediate covering $X_\chi(N)\longrightarrow X_0(N)$
of the degeneracy morphism $\pi$ associated to the subgroup $\Omega$.
In fact $X_\chi(N)=\Gamma_\chi(N)\backslash\mathcal{H}^*$,
where 
\begin{equation*}
\Gamma_\chi(N)=
\left\{\mtwo{\alpha}{\beta}{\gamma}{\delta}\in\Gamma_0(N)\,\colon\,
\chi(\delta)= 1\right\}.
\end{equation*}
The curve $X_\chi(N)$ has $4$ cups,
namely the cusps $\infty_1$ and $\infty_2$,
above the cusp $\infty$ and cusps $0_1$ and $0_2$,
above the cusp $0$;
the Galois group $G(X_\chi(N)/X_\chi(N))$
acts transitively on the set $\{\infty_1,\infty_2\}$
and on the set $\{0_1,0_2\}$.

%
\section{Variant of Ligozat's construction}\label{sec:lig}
Recall that for odd primes $p$ and $q$ such that $p\not=q$ Gau\ss\
proved that 
\begin{equation*}
\left(\frac{q}{p}\right)=(-1)^{S(q,p)}
\end{equation*}
where $(\frac{\cdot}{\cdot})$ is the Legendre symbol,
\begin{equation*}
S(q,p)=\sum_{r=1}^{\frac{p-1}{2}}\left\lfloor r\frac{q}{p}\right\rfloor
\end{equation*}
and $\lfloor x \rfloor$ is the usual floor function.
(Cf. p. 78 of Hardy and Wright~\cite{hardy:intro}.)
This result may be generalised as follows.  

\begin{lem}\label{lem:sun}
Suppose $n$ is an odd integer.
If $\delta$ is any integer not divisible by $n$ then
\begin{equation*}
\left(\frac{\delta}{n}\right)=
(-1)^{S(\delta,n)}(-1)^{\frac{n^2 - 1}{8}(\delta+1)}
\end{equation*}
where $(\frac{\cdot}{\cdot})$ is the Jacobi symbol.
\end{lem}

\begin{proof}
On the one hand by the work of Jenkins~\cite{jenkins:reciprocity}
we know that $(\frac{\delta}{n})=(-1)^{\nu(\delta,n)}$,
where $\nu(\delta,n)$ is the number of elements of the set
\begin{equation*}
\left\{k\in\ZZ_{>0}\colon k < \frac{n}{2} <(k\delta\pmod{n})\right\}.
\end{equation*}
On the other hand Lemma 3.1 of Zhi-Wei Sun~\cite{sun:binomial} implies
the congruence
\begin{equation*}
 \nu\equiv S(\delta, n) + \frac{n^2 -1}{8}(\delta +1)\pmod{2}
\end{equation*}
and the lemma follows.
\end{proof}

\begin{thm}[Ogg, Ligozat]\label{thm:modunit}
There is a holomorphic function on $\hyper$ such that
\begin{equation}\label{eqn:etaiota}
f_\chi(\mu\tau)=\chi(\delta)f_\chi(\tau)^{\chi(\delta)},
\end{equation}
where $\mu=\abcd{\alpha}{\beta}{\gamma}{\delta}\in\Gamma_0(N)$.
\end{thm}

\begin{proof}
Following Kubert and Lang~\cite{lang:kl02} we
define the \textit{Siegel function}
\begin{equation*}
\frak{g}_a(\tau)=\frak{t}_a\vtwo{\tau}{1}\Delta^\frac{1}{12}(\tau)
\end{equation*}
where $\frak{g}_a$ is the Klein form associated to a
non-zero element $a=(a_1,a_2)$
of $(\frac{1}{N}\ZZ/\ZZ)\times(\frac{1}{N}\ZZ/\ZZ)$,
and $\Delta^\frac{1}{12}(\tau)$ is
the principal part of the $12$-\textit{th} root of the cusp form
\begin{equation*}
\Delta(\tau)=q\prod_{n=1}^\infty(1-q^n)^{24}\qquad(q=e^{2\pi i\tau})
\end{equation*}
of weight $12$.
It is well-known that
\begin{equation*}
\frak{g}_a(\mu\tau)=\xi(\mu)\frak{g}_{a\mu}(\tau)
\end{equation*}
where $\xi(\mu)$ is a the $12$-\textit{th} root of unity
such that $\xi\abcd{0}{1}{-1}{0}=i$
and $\xi\abcd{1}{1}{0}{1}=\zeta_{12}$.
(See Serre's book~\cite{serre:arithmetic}.)
Now for each integer $r$ such that gcd$(r,N)=1$ we define
\begin{equation*}
f_r(\tau)=\frak{g}_{(0,\frac{r}{N})}(\tau)
\end{equation*}
So the function of Ogg and Ligozat
(cf. Mazur~\cite{mazur:eisen}, p. 107)
may be defined by the product
\begin{equation*}
f_\chi(\tau)=
\prod_{r=1}^{\frac{N-1}{2}} f_r(\tau)^{\chi(r)}.
\end{equation*}
By the work of Kubert and Lang~\cite{lang:munits} we know that 
\begin{equation*}
f_r(\mu\tau)=\xi(\mu)
e^{\pi i(-\frac{r^2}{N^2}\gamma\delta + 
\frac{r}{N}\gamma + \lfloor \frac{r}{N}\delta \rfloor)}f_{\{r\delta\}}(\tau)
\end{equation*}
So we have
\begin{equation*}
f_\chi(\mu\tau)=
\prod_{r=1}^{\frac{N-1}{2}}f_r(\mu\tau)^{\chi(r)}=
e^{\pi i(S_1 + S_2 + S_3)}f_\chi(\tau)^{\chi(\delta)},
\end{equation*}
where
\begin{equation*}
S_1 = -\frac{\gamma\delta}{N^2}\sum_{r=1}^{\frac{N-1}{2}}\chi(r)r^2
\equiv\delta\frac{\gamma}{N}\cdot\frac{N^2-1}{24}
\equiv\delta\gamma\cdot\frac{N^2-1}{8}
\pmod{2},
\end{equation*}
\begin{equation*}
S_2 = \frac{\gamma}{N}\sum_{r=1}^{\frac{N-1}{2}}\chi(r)r
\equiv\frac{\gamma}{N}\cdot\frac{N^2-1}{8}
\equiv\gamma\cdot\frac{N^2-1}{8}\pmod{2},
\end{equation*}
\begin{equation*}
S_3 = \sum_{r=1}^{\frac{N-1}{2}}
\chi(r)\left\lfloor \frac{r}{N}\delta \right\rfloor
\equiv (\delta + 1)\cdot\frac{N^2-1}{8} + k_{N,\delta} \pmod{2},
\end{equation*}
and $k_{\delta,N}=1$,
if $(\frac{\delta}{N})=-1$ and $k_{\delta,N}=0$ otherwise.
(The latter congruence follows from Lemma~\ref{lem:sun}.)
Therefore
\begin{equation*}
f_\chi(\mu\tau)=\chi(\delta)(-1)^{(\delta+1)(\gamma+1)\frac{N^2-1}{8}}
f_\chi(\tau)^{\chi(\delta)}
\end{equation*}
The equation $\alpha\delta - \beta\gamma=1$ implies
that $\gamma$ and $\delta$ may not have the same parity. 
So the theorem follows.
\end{proof}

%
\section{Infinite products and Fourier expansions}
The \textit{generalised Bernoulli numbers} $B_{n,\chi}$
attached to the non-trivial primitive character $\chi$
of conductor $N$ may be defined by the formal power series
\begin{equation*}
\sum_{r=1}^N \chi(r)\frac{Xe^{rX}}{e^{NX}-1}=
\sum_{n=0}^\infty B_{n,\chi}\frac{X^n}{n!}
\end{equation*}

\begin{lem}\label{lem:cyclotomic}
The quotient of polynomials
\begin{equation*}
\Psi(X)=\prod_{r=1}^{N-1}(1-X\zeta_N^r)^{\chi(r)}
\end{equation*}
has a formal power series expansion of the form
\begin{equation*}
\Psi(X)= 1 - \sqrt{N}X + \dots \in k[[X]],
\end{equation*}
where $k=\QQ(\sqrt{N})$.
\end{lem}

\begin{proof}
Let $H$ denote the (unique) subgroup of index $2$ of
the Galois group $G(\QQ(\zeta_N)/\QQ)$ and note that our
assumption $N\equiv 1\pmod{4}$ implies that
the fixed field of $H$ is the real quadratic field $k$.
Now let $\sigma$ be an element of $H$ and
let $\rho$ be the element of $(\ZZ/N\ZZ)^\times$ with the property 
$\sigma(\zeta_N)=\zeta_N^\rho$.
Clearly
\begin{equation*}
\sigma(\Psi(X))=\sigma\left(\prod_{r=1}^{N-1}(1-X\zeta_N^r)^{\chi(r)}\right)=
\prod_{r=1}^{N-1}(1-X\zeta_N^{\rho r})^{\chi(r)}=\Psi(X).
\end{equation*}
So $\Psi(X)$ is a quotient of two polynomials defined over
$\QQ(\sqrt{N})$, and thus
the formal power series of $\Psi(X)$ has coefficients in $\QQ(\sqrt{N})$.
Finally note that
\begin{displaymath}
(1-X\zeta_N^r)^{\chi(r)}=\left\{
\begin{array}{ll}
1 + \zeta_N^rX+\zeta_N^{2r}X^2+\dots,&\textit{if $\quad\chi(r)=-1$,}\\
1 - \zeta_N^rX,&\textit{if $\quad\chi(r) = +1$.}\\
\end{array}
\right.
\end{displaymath}
So clearly 
\begin{equation*}
\Psi(X) = 1 -\left(\sum_{r=1}^{N-1}\chi(r)\zeta^r\right)X + \dots=
1 -\sqrt{N}X+\dots
\end{equation*}
where the last equality follows from
a well-known result of Gau\ss\
and the fact that $N\equiv 1\pmod{4}$.
\end{proof}

As before put $\breve{f}_\chi(\tau)=w_Nf_\chi(\tau)=f_\chi(w_N(\tau))$,
where $w_N$ is the Fricke involution $w_N(\tau)=-\frac{1}{N\tau}$.

\begin{thm}\label{thm:infnprod}
Assume for simplicity that $N>5$.
We have the following $q$-expansions.
\begin{description}
\item[A] If $k=\QQ(\sqrt{N})$ and
$h(N)$ is the class number of $k$, then
$$f_\chi(\tau)=u^{h(N)}(1-\sqrt{N}q+\dots)\in k[[q]],$$
where $u$ is a fundamental unit of $k$ and $h(N)$ is the class number of $k$.
\item[B] If as above $B_{2,\chi}$ denotes the second Bernoulli number
attached to $\chi$ then
$$\breve{f}_\chi(\tau)=
q^{\frac{1}{2}B_{2,\chi}}\prod_{n=1}^\infty(1-q^n)^{\chi(n)}\in \ZZ[[q]]$$
\end{description}
\end{thm}

\begin{proof}
By well-known results on Siegel functions
(see Kubert and Lang~\cite{lang:munits})
and Lemma~\ref{lem:cyclotomic} we have
\begin{equation*}
f_\chi(\tau)=
\prod_{r=1}^{\frac{N-1}{2}}
\left(
\frac{-1}{2\pi i}q^{\frac{1}{12}}
\zeta_{2N}^{-r}(1 - \zeta_N^r)
\prod_{n=1}^\infty(1 - q^n\zeta_N^r)(1 - q^n\zeta_N^{-r})
\right)^{\chi(r)}=
\end{equation*}
\begin{equation*}
\lambda\prod_{r=1}^{\frac{N-1}{2}}\prod_{n=1}^\infty
(1 - q^n\zeta_N^r)^{\chi(r)}(1 - q^n\zeta_N^{-r})^{\chi(-r)}=
\lambda\prod_{n=1}^\infty\prod_{r=1}^{N-1}
(1 - q^n\zeta_N^r)^{\chi(r)}=
\end{equation*}
\begin{equation*}
\lambda\prod_{n=1}^\infty\Psi(q^n),
\end{equation*}
where
\begin{equation*}
\lambda=
\prod_{r=1}^{\frac{N-1}{2}}
\zeta_{2N}^{-\chi(r)r}(1-\zeta_N^r)^{\chi(r)}.
\end{equation*}
So it remains to show that $\lambda=u^{h(N)}$.
Note
\begin{equation*}
\lambda^\kappa=
\left(\prod_{r=1}^{\frac{N-1}{2}}
(\zeta_N^{-\frac{r}{2}}-\zeta_N^{\frac{r}{2}})^{\chi(r)}\right)^\kappa=
(-1)^{\frac{N-1}{2}}\lambda=\lambda
\end{equation*}
where $\kappa$ denotes complex conjugation.
So $\lambda$ is real.
Moreover
\begin{equation*}
\lambda^2=\lambda\lambda^\kappa=
\prod_{r=1}^{\frac{N-1}{2}}(1-\zeta_N^r)^{\chi(r)}
\prod_{r=1}^{\frac{N-1}{2}}(1-\zeta_N^{-r})^{\chi(-r)} =
\prod_{r=1}^{N-1}(1-\zeta_N^r)^{\chi(r)}.
\end{equation*}
By the work of Tate~\cite{tate:stark}
(cf. Darmon~\cite{darmon:stark}) we know
\begin{equation*}
\prod_{r=1}^{N-1}(1-\zeta_N^r)^{\chi(r)}=u^{2h(N)},
\end{equation*}
where $u$ is a fundamental unit of the real quadratic field $k=\QQ(\sqrt{N})$.
Therefore $\lambda=u^{\pm h(N)}$ and \textbf{(A)} follows.
Now consider $N\equiv 1\pmod{4}$ implies
\begin{equation*}
\frac{N}{2}\sum_{r=1}^{\frac{N-1}{2}}
\chi(r)B_2\left(\frac{r}{N}\right)=\frac{1}{2}B_{\chi,2},
\end{equation*}
where $B_2(X)$ is the second Bernoulli polynomial introduced earlier.
Moreover,
\begin{equation*}
f_\chi(\tau)=
\prod_{r=1}^{\frac{N-1}{2}}
\left(
\frac{-1}{2\pi}q^{\frac{1}{2}B_2(\frac{r}{N})N}
(1-q^r)\prod_{n=1}^\infty (1 - q^{Nn + r})(1 - q^{Nn - r})
\right)^{\chi(r)}=
\end{equation*}
\begin{equation*}
q^{\frac{1}{2}B_{2,\chi}}
\prod_{r=1}^{\frac{N-1}{2}}(1-q^r)
\prod_{r=1}^{\frac{N-1}{2}}
\prod_{n=1}^\infty (1-q^{Nn + r})^{\chi(r)}(1-q^{Nn - r})^{\chi(-r)}=
\end{equation*}
\begin{equation*}
q^{\frac{1}{2}B_{2,\chi}}\prod_{n=1}^\infty(1-q^n)^{\chi(n)},
\end{equation*}
and \textbf{(B)} follows.
\end{proof}

From Theorem~\ref{thm:infnprod} and Theorem~\ref{thm:modunit}
it follows that $f_\chi$ defines a rational function on $X_\chi(N)$,
provided that $\frac{1}{2}B_{n,\chi}$ is integral.
In fact $\frac{1}{2}B_{n,\chi}$ is integral for all primes $N$
except when $N=5$,
where $\frac{1}{2}B_{n,\chi}=\frac{1}{5}$.
In this special case we will find it convenient to switch to the notation 
\begin{equation*}
f_\chi(\tau)=
\left(\prod_{r=1}^{\frac{N-1}{2}} f_r(\tau)^{\chi(r)}\right)^5.
\end{equation*}

%
\section{On some identities of Ramanujan}\label{sec:ram}
Suppose $N\equiv 1\pmod{4}$ prime.
Let $\breve{g}_\chi = \breve{f}_\chi - \frac{1}{\breve{f}_\chi}$.
It follows from Theorem~\ref{thm:modunit} that $\breve{g}_\chi$
defines a rational function on $X_0(N)$.
Also let
\begin{equation*}
t(\tau)=q^{\frac{1-N}{m}}
\prod_{\genfrac{}{}{0pt}{}{n>0}{n\nmid N}}(1-q^n)^{\frac{24}{m}}
\end{equation*}
where $m=\textrm{gcd}(N-1,12)$.
It is well-known that $t(\tau)$ defines a rational function on
the modular curve $X_0(N)$.
Certain identities due to Ramanujan
may be succinctly expressed in terms of $\breve{g}_\chi$ and $t(\tau)$.
Indeed,
suppose $N=5$ and consider
Entry 11(iii) of Chapter 19 of Berndt's book~\cite{berndt:notebook}
may be expressed as
\begin{equation*}
((11+t(\tau))^2+1)^{1\over2}-(11+t(\tau))=\breve{f}_\chi(\tau),
\end{equation*}
and thus
\begin{equation*}
\breve{f}_\chi^{2}(\tau) +(11 + t(\tau))\breve{f}_\chi(\tau) -1 = 0.\footnote{
The function $\breve{f}_\chi(\tau)$ also appears in connection with
the Rogers-Ramanujan's continued fractions,
e.g. $1+\frac{e^{-2\pi}}{1+}\,\frac{e^{-4\pi}}{1+}\,\frac{e^{-6\pi}}{1+\dots}=
\left(
\sqrt{\frac{5+\sqrt{5}}{2}} -\frac{1+\sqrt{5}}{2}
\right)e^{\frac{2\pi}{5}}$
in $\S 19.15$ (p. 294-295) of Hardy and Wright~\cite{hardy:intro}.}
\end{equation*}
In other words
\begin{equation}\label{eqn:ram05}
\breve{g}_\chi(\tau) = 11 + t(\tau).
\end{equation}
Now suppose $N=13$.
Then $\mu_2\mu_3\mu_4 = \breve{f}_\chi$
and $\mu_1\mu_5\mu_6 = \frac{1}{\breve{f}_\chi}$,
where
\begin{align*}
\mu_1&=q^{-\frac{7}{13}}\frac{f(-q^4,-q^9)}{f(-q^2,-q^{11})},&\mu_2=q^{-\frac{6}{13}}\frac{f(-q^6,-q^7)}{f(-q^3,-q^{10})},\\
\mu_3&=q^{-\frac{5}{13}}\frac{f(-q^2,-q^{11})}{f(-q,-q^{12})},&\mu_4=q^{-\frac{2}{13}}\frac{f(-q^5,-q^8)}{f(-q^4,-q^9)},\\
\mu_5&=q^{ \frac{5}{13}}\frac{f(-q^3,-q^{10})}{f(-q^5,-q^8)},&\mu_6=q^{\frac{15}{13}}\frac{f(-q,-q^{12})}{f(-q^6,-q^7)},\\
\end{align*}
and $f(a,b)$ is Ramanujan's two-variable theta function
\begin{equation*}
f(a,b)=\prod_{n=0}^\infty(1+a^{n+1}b^n)(1+a^nb^{n+1})(1-a^{n+1}b^{n+1}).
\end{equation*}
So Entry 8 of Chapter 20 of Berndt's book~\cite{berndt:notebook}
\begin{equation*}
t + 3= \mu_2\mu_3\mu_4 - \mu_1\mu_5\mu_6
\end{equation*}
yields
\begin{equation*}
\breve{f}_\chi^2(\tau) + (t(\tau) + 3)\breve{f}_\chi(\tau) - 1 = 0.
\end{equation*}
Hence
\begin{equation}\label{eqn:ram13}
\breve{g}_\chi(\tau) = t(\tau) + 3.
\end{equation}

From the point of view of the function theory of the curve $X_0(N)$
Equation~\ref{eqn:ram05} and Equation~\ref{eqn:ram13} are essentially obvious.
Indeed, 
consider the well-known fact that the divisor of poles of $t(\tau)$
is concentrated at the cusp $\infty$ with degree $v_N=\frac{1-N}{m}$,
and that the divisor of poles of $g_\chi(\tau)$
is concentrated also at the cusp $\infty$,
but with degree $v_\chi=\frac{1}{2}B_{2,\chi}$.
Assume $N=5,13$.
Thus both,
the modular unit $t(\tau)$ and the rational function $g_\chi(\tau)$
have valence one.
So there must be a polynomial $P(T)\in\ZZ[T]$ of degree one
such that $g_\chi(\tau)=P(t(\tau))$,
which may be explicitly obtained using the principal part
of the $q$-expansion of each of these functions.

An equation of the form $g_\chi(\tau)=P(t(\tau))$,
for some polynomial $P(T)\in\QQ[T]$,
for prime level $N\equiv 1\pmod{4}$ greater than $13$ seems unlikely.
Indeed,
experimental evidence suggests that for such levels $N>13$
the necessary divisibility condition $v_N\mid v_\chi$ holds
are exactly the ones given by the table:
\begin{displaymath}
\begin{array}{lrr}
N&  v_N&  v_\chi\\
\vspace{-2ex}\\
109& 9&  27\\
197&  49&  49\\
2617&  218&  4796\\
3709&  309&  5253\\
\end{array}
\end{displaymath}
By subtracting from $g_\chi(\tau)$ the relevant powers of $t(\tau)$,
for each level $N$ in the above table
it is easy to see that a relation of the form $g_\chi(\tau)=P(t(\tau))$
is not possible.

For $X_0(N)$ of genus $g_N=1$ it is not difficult compute explicit
expressions of the form $g_\chi(\tau)=P(X(\tau),X(\tau))$,
where $P(X,Y)\in\QQ[X,Y]$,
and $X$ and $Y$ are given generators of
the function field of $X_0(N)$ over $\QQ$.
Indeed,
let $X$ and $Y$ be the coordinate functions of the
modular parametrisation of $X_0(N)$,
regarded as an elliptic curve.
Then by computing a few terms of the $q$-expansion of $g_\chi(\tau)$
(see Theorem~\ref{thm:infnprod})
it is easy to find an integer $a(r,s)$,
and positive integers $r$ and $s$ such that
the rational function  $g_\chi(\tau)-a(r,s)X^rY^s$
has a pole at $\infty$ of degree less than that of $g_\chi(\tau)$ at $\infty$.
Proceeding recursively,
we may arrive at a function with no poles that vanishes at $\infty$,
thus the zero function on $X_0(N)$.
This yields the  desired expression $g_\chi=P(X,Y)$.

\begin{example}
Suppose $N = 17$.
Using Cremona's Tables~\cite{cremona:onlinetables} we may see that
the curve $X_0(N)$ of genus $g_N = 1$ is isomorphic to curve {\bf 17A1}.
Thus it has a global minimal Weierstra\ss\ model
\begin{equation*}
Y^2 + XY + Y = X^3 - X^2 - X - 14,
\end{equation*}
and a modular parametrisation
\begin{align*}
X(\tau)&= \frac{1}{q^2} + \frac{1}{q} + 1 + O(q),\\ 
Y(\tau)&=-\frac{1}{q^3} - \frac{2}{q^2} - \frac{2}{q} - 3 + O(q).
\end{align*}
Since $\breve{g}_\chi$ has exactly one pole,
and $\breve{g}_\chi = -\frac{1}{q^2} - \frac{1}{q} - 2 + O(q)$,
it follows that
\begin{equation*}
\breve{g}_\chi(\tau) = - (X(\tau) + 1).
\end{equation*}
\end{example}

For $X_0(N)$ of genus $g_N=2$ the situation is nearly as simple.

\begin{example} 
Suppose $N=29$.
The modular curve $X_0(N)$ is
an hyperelliptic curve of genus $g_N=2$ with a model
\begin{equation*} 
Y^2 = X^6 + 2X^5 - 17X^4 - 66X^3 - 83X^2 - 32X - 4,
\end{equation*} 
where $X= \sqrt{2}\cdot\frac{f^\sigma + f}{f^\sigma - f}$
and $Y= q\left(\frac{d}{dq}X\right)\frac{1}{f^\sigma - f}$,
and $f$ is the (unique) newform of level $N$ weight $2$ and trivial character.
With the help of the fact that $w_NX=X$ and $w_NY=-X$ it may be found that
\begin{equation*} 
\breve{g}_\chi(\tau) = -\frac{1}{2}(Y + X^3 + X^2 - 9X) + 7.
\end{equation*}
\end{example}

\begin{example}
Suppose $N=37$.
It is well-known that the curve $X_0(N)$ has genus $g_N=2$ and a model
\begin{equation*} 
Y^2 = X^6 + 14X^5 + 35X^4 + 48X^3 + 35X^2 + 14X + 1
\end{equation*}
where $X = \frac{f_2 + f_1}{f_2 - f_1}$
and $Y= 2q\left(\frac{d}{dq}X\right)\frac{1}{f_1 - f_2}$,
and $f_1$ (resp. $f_2$) is the rational newform for $\Gamma_0(N)$
such that the Fourier coefficient $a_{f_1}(N) = -1$ (resp. $a_{f_2}(N) = 1$).
Then $\breve{g}_\chi(\tau)$ is $-\frac{1}{2}$ times
%
%
%
\begin{equation*}
X^5 +16X^4 +67X^3 +YX^2 +87X^2 +9YX +62X +11Y +13
\end{equation*}
\end{example}

%
\section{The zeros of $h_\chi$ for small $N$}\label{sec:hchi}
As above fix a prime $N\equiv 1\pmod{4}$.
Define
\begin{equation*}
h_\chi = g_\chi + \breve{g_\chi} = f_\chi - \frac{1}{f_\chi} +
\breve{f}_\chi - \frac{1}{\breve{f}_\chi}.
\end{equation*}
A consequence of Theorem~\ref{thm:modunit} is that $h_\chi$ may be
regarded as a rational function on $X_0^+(N)$.
Using again basic results from the Theory of Algebraic Curves
it is an easy matter to obtain
equations of the form $h_\chi(\tau)=P(X(\tau),X(\tau))$,
where $P(X,Y)\in\QQ[X, Y]$,
and $X$ and $Y$ are generators of
the function field of $X_0^+(N)$ over $\QQ$.
The polynomial $P(X,Y)$ may be obtained in a completely analogous way
as we did for the rational function $\breve{g}_\chi(\tau)$ on $X_0(N)$,
but now using the formula
\begin{equation*}
h_\chi(\tau) = -T_{k/\QQ}(u^{h(N)}) -\frac{1}{\breve{f_\chi}} + O(q),
\end{equation*}
where $u$ is a fundamental unit of $k$,
and $h(N)$ is the class number of $k$,
and $T_{k/\QQ}$ is the trace from the real quadratic field $k=\QQ(\sqrt{N})$
to the field of the rational numbers $\QQ$.
(See Theorem~\ref{thm:infnprod}.)
Once we know $P(X,Y)$ we may (at least in principle) obtain the field of
definition of each of the zeros of $h_\chi(\tau)$ by
decomposing into irreducible factors over $\QQ$
a generator $g_X$ of the kernel of the natural map
\begin{equation*}
\QQ[X, Y]/(W,P)\longrightarrow\QQ[X, Y]/(Y)
\end{equation*}
and a generator $g_Y$ of the kernel of the natural map
\begin{equation*}
\QQ[X, Y]/(W,P)\longrightarrow\QQ[X, Y]/(X),
\end{equation*}
provided $X_0^+(N)$ is embeddable in the projective plane $\PP^2$
and that we have an explicit equation $W(X,Y)=0$ for $X_0^+(N)$.

\begin{example}
Let $E$ be elliptic curve \textbf{37A1}.
It is the elliptic curve of conductor $N=37$
with a minus sign in the functional equation of
its $\Lambda$-function.
(See Cremona's Tables~\cite{cremona:onlinetables}.)
So $X_0^+(N)$ may be identified with $E$ via
the modular parametrisation $\tau\mapsto(X(\tau),Y(\tau))$,
where\footnote{
See Zagier's paper~\cite{zagier:modular} for
an excellent account on the construction of this modular parametrisation.}
\begin{equation*}
\begin{split}
X(\tau)&= \frac{1}{q^2} + \frac{2}{q} + 5 + O(q) \\
Y(\tau)&=-\frac{1}{q^3} - \frac{3}{q^2} - \frac{9}{q} - 21 + O(q)
\end{split}
\end{equation*}
and $q=e^{2\pi i\tau}$.
Now consider that the modular unit $\breve{f}_\chi$ has divisor
\begin{equation*}
(\breve{f}_\chi)=\frac{1}{2}B_{\chi,2}(\infty_1-\infty_2)=5(\infty_1-\infty_2)
\end{equation*}
Thus the poles of the function $h_\chi(\tau)$ are
concentrated at the only cusp of $X_0^+(N)$,
with multiplicity $5$.
So the Riemann-Roch theorem implies that $h_\chi$
may be expressed as a (unique) linear combination of
the $5$ functions $1,X,Y,X^2$ and $XY$.
So a simple exercise in linear algebra involving
the principal parts of the $q$-expansion
\begin{equation*}
h_\chi(\tau)=
-\frac{1}{q^5} -\frac{1}{q^4} -\frac{1}{q^2} - \frac{2}{q} - 12 + O(q),
\end{equation*}
(see Theorem~\ref{thm:infnprod})
and of the $q$-expansions of the above $5$ functions yield
\begin{equation*}
h_\chi= 1 + 6X + 4Y - 4X^2 - XY.
\end{equation*}
Now with the help of Grayson's \textsc{Macaulay2}~\cite{grayson:m2},
the latter equation together with the Weierstra\ss\ equation of $E$:
\begin{equation*}
Y^2 + Y = X^3 - X,
\end{equation*}
imply that
the $X$-coordinates of the zeros of $h_\chi$ are exactly the roots of
\begin{equation*}
g_X(T)=T^5 - 24T^4 + 67T^3 - 42T^2 - 5T + 3 \in \QQ[T].
\end{equation*}
The polynomial $p_x(T)$ factors into irreducibles (over $\QQ$) as
\begin{equation*} 
g_X(T)=(T-1)(T^4 - 23T^3 + 44T^2 + 2T - 3).
\end{equation*}
Similarly,
the $Y$-coordinates of the zeros of $h_\chi$ are exactly the roots of
\begin{equation*}
g_Y(T)=T^5  + 95T^4  - 86T^3  - 279T^2  - 72T + 27,
\end{equation*}
which factors into irreducibles (over $\QQ$) as
\begin{equation*}
g_Y(T)=(T + 1)(T^4 + 94T^3 - 180T^2 - 99T + 27).
\end{equation*} 
\end{example}

The computations for $X_0^+(N)$ of genus $g_N^+=2$
are nearly as simple as for $g_N^+=1$,
so we do not discuss the details of these computations here any further.
For the $N=53$,  $61$, $89$, and $101$ (i.e. the rest of the $g_N^+=1$ cases)
and for $N=73$ (the only $g_N^+=2$ case under our hypothesis)
we found exactly one nontrivial Galois orbit of roots of $h_\chi(\tau)$.
The minimal polynomial $p_N(X)\in\ZZ[T]$
(normalised so that its coefficients are coprime)
of a generator of the
field of definition of these points is given as follows.
\begin{displaymath}
p_{53}(T)=T^6 - 20T^5 + 95T^4 - 156T^3 + 145T^2 - 174T - 44.
\end{displaymath}
\begin{displaymath}
\begin{split}
p_{61}(T)=&T^{10} +72T^9 +1000T^8 -2327T^7 -2810T^6 +3994T^5\\
       & +1905T^4 -2283T^3 -221T^2 + 404T -60
\end{split}
\end{displaymath}
\begin{displaymath}
\begin{split}
p_{89}(T)=&T^{26} -22T^{25}+108T^{24}+489T^{23}-3164T^{22}-9330T^{21}\\
       &+30025T^{20} +120140T^{19}-4651T^{18}-483581T^{17}\\
       &-576269T^{16}+246025T^{15}+882959T^{14}+596485T^{13}\\
       &+263186T^{12}+59943T^{11}-289362T^{10}-263968T^9\\
       &-43576T^8-51782T^7-44804T^6+3630T^5\\
       &+7476T^4 -6260T^3-4128T^2+3680T-352
\end{split}
\end{displaymath}
\begin{displaymath}
\begin{split}
p_{101}(T)=&T^{18} - 19T^{17} + 135T^{16} - 434T^{15} + 548T^{14} + 145T^{13}\\
       &- 1028T^{12} + 1631T^{11} - 2464T^{10} + 1016T^9 + 4005T^8\\
       &- 6040T^7 - 1811T^6 + 5457T^5 - 55T^4 - 3417T^3\\
       &- 195T^2 + 1694T + 676
\end{split}
\end{displaymath}
\begin{displaymath}
\begin{split}
p_{73}(T)=&8T^{22} - 84T^{21} - 874T^{20} + 142T^{19} + 15945T^{18} + 26187T^{17}\\
       &- 98676T^{16} - 300010T^{15} + 117375T^{14} + 1211979T^{13}\\
       &+ 802441T^{12} - 1804645T^{11} - 2343277T^{10} + 714633T^9\\
       &+ 1965510T^8 - 93748T^7 - 882954T^6 + 62476T^5 + 225574T^4\\
       &- 47106T^3 - 22652T^2 + 9508T - 984
\end{split}
\end{displaymath}

Note that the degree of each of the above polynomials
is one less than the degree of the divisor of zeros
of $h_\chi(\tau)$ for $N=37$, $53$, $61$, and $101$.
These are exactly the prime levels such that $N\equiv 5\pmod{8}$
considered just above.
This is explained by the existence of a ``trivial''
rational zero of $h_\chi(\tau)$.

\begin{lem}\label{lem:zero}
If $N\equiv 5\pmod{8}$ then $h_\chi(\tau)$
vanishes at the common image in $X_0^+(N)$ of
the two Heegner points of discriminant $D=-4$ on the curve $X_0(N)$.
\end{lem}

\begin{proof}
The congruence $N\equiv 5\pmod{8}$ implies that
the non-trivial element of the
Galois group $G(X_\chi(N)/X_0(N))$ may be represented by the matrix
\begin{equation*}
\mu=\mtwo{\rho}{1}{-(\rho^2+1)}{-\rho}\in\Gamma_0(N),
\end{equation*}
where $\rho$ is any integer such that $\rho^2\equiv -1\pmod{N}$.
Note that $\mu$ is an elliptic matrix of order $2$ that fixes the 
Heegner points $\tau_{-4,\pm r}=-\frac{1}{i \pm r}\in\hyper$,
of discriminant $D=-4$ of $X_0(N)$.
Equation~\ref{eqn:etaiota} yields $\breve{f}_\chi(\tau_{-4, \pm r})^2=-1$.
But  $\breve{f}_\chi(\tau)^\kappa=\breve{f}_\chi(\tau^\kappa)$,
where $\kappa$ denotes complex conjugation.
Thus $\breve{f}(\tau_{-4,r})=i$ implies
\begin{equation*}
\breve{g}_\chi(\tau_{-4,-r})=
\breve{f}_\chi(\tau_{-4,-r})-\frac{1}{\breve{f}_\chi(\tau_{-4,-r})}=-2i.
\end{equation*}
Similarly,
if $\breve{f}(\tau_{-4,r})=-i$ then
\begin{equation*}
\breve{g}_\chi(\tau_{-4,-r})=
\breve{f}_\chi(\tau_{-4,-r})-\frac{1}{\breve{f}_\chi(\tau_{-4,-r})}=2i.
\end{equation*} 
But $w_N(\tau_{-4,r})=\tau_{-4,-r}$.
Therefore $\breve{g}_\chi(\tau_{-4,\pm r})=-g_\chi(\tau_{-4,\pm r})$ and
the lemma follows.
\end{proof}

%
\section{Final remarks}\label{sec:conj}
Some computations we performed on the
above polynomials $p_N(T)$ are consistent with the following.

\begin{conjecture}
The field of definition
each of the zeros of $h_\chi$ is either $\QQ$ or
an extension $F$ of the real quadratic field $k=\QQ(\sqrt{N})$
of degree $n=\frac{1}{2}B_{2,\chi}$,
if $N\equiv 1\pmod{8}$ and $n=\frac{1}{2}B_{2,\chi}-1$,
if $N\equiv 5\pmod{8}$,
where $B_{m,\chi}$ is the $m$-th generalised Bernoulli number
attached to the character $\chi$.
\end{conjecture}

In other words we conjecture that the Galois group is in
a sense as ``large as possible''.
Numerical evidence suggests that
the above conjecture is true also for the genus $g_N^+=0$ case,
i.e. for $N=17$, $29$, and $N=41$.
It is hoped that a more extensive numerical evidence from
higher genus examples, e.g. $g_N^+=3$
may perhaps shed some light on the nature of phenomena.

Finally,
perhaps it is worth investigating if this potential connexion with
Singmaster's Problem~\cite{singmaster:problem}
might suggest a yet to be discovered,
interesting combinatorial aspect of the Galois group $G=G(F^\textit{nrm}/\QQ)$
of the normal closure $F^\textit{nrm}$ of $F$ over the rationals $\QQ$.


\bibliographystyle{amsplain}
\bibliography{biblio}
\printindex
\end{document}